\documentclass[12pt]{article}

\usepackage{amsmath,amssymb,amsthm,amscd,a4wide,upref}
\usepackage[enableskew]{youngtab}
\usepackage{hyperref}
\hypersetup{colorlinks,citecolor=blue,filecolor=black,linkcolor=blue,urlcolor=blue}
\usepackage{enumerate}
\usepackage{mathtools}
\usepackage{enumitem}
\usepackage{lmodern}     
\usepackage[T1]{fontenc}

\begin{document}

\newcommand{\End}{{\rm{End}\ts}}
\newcommand{\Hom}{{\rm{Hom}}}
\newcommand{\Mat}{{\rm{Mat}}}
\newcommand{\ch}{{\rm{ch}\ts}}
\newcommand{\sh}{{\rm{sh}}}
\newcommand{\chara}{{\rm{char}\ts}}
\newcommand{\diag}{ {\rm diag}}
\newcommand{\non}{\nonumber}
\newcommand{\wt}{\widetilde}
\newcommand{\wh}{\widehat}
\newcommand{\ot}{\otimes}
\newcommand{\la}{\lambda}
\newcommand{\La}{\Lambda}
\newcommand{\De}{\Delta}
\newcommand{\al}{\alpha}
\newcommand{\be}{\beta}
\newcommand{\ga}{\gamma}
\newcommand{\Ga}{\Gamma}
\newcommand{\ep}{\epsilon}
\newcommand{\ka}{\kappa}
\newcommand{\vk}{\varkappa}
\newcommand{\si}{\sigma}
\newcommand{\vp}{\varphi}
\newcommand{\de}{\delta}
\newcommand{\ze}{\zeta}
\newcommand{\om}{\omega}
\newcommand{\ee}{\epsilon^{}}
\newcommand{\su}{s^{}}
\newcommand{\hra}{\hookrightarrow}
\newcommand{\ve}{\varepsilon}
\newcommand{\ts}{\,}
\newcommand{\vac}{\mathbf{1}}
\newcommand{\vacr}{|\tss 0\rangle}
\newcommand{\vacl}{\langle 0\tss |}
\newcommand{\di}{\partial}
\newcommand{\qin}{q^{-1}}
\newcommand{\tss}{\hspace{1pt}}
\newcommand{\Sr}{ {\rm S}}
\newcommand{\U}{ {\rm U}}
\newcommand{\BL}{ {\overline L}}
\newcommand{\BE}{ {\overline E}}
\newcommand{\BP}{ {\overline P}}
\newcommand{\AAb}{\mathbb{A}\tss}
\newcommand{\CC}{\mathbb{C}\tss}
\newcommand{\KK}{\mathbb{K}\tss}
\newcommand{\QQ}{\mathbb{Q}\tss}
\newcommand{\SSb}{\mathbb{S}\tss}
\newcommand{\ZZ}{\mathbb{Z}\tss}
\newcommand{\X}{ {\rm X}}
\newcommand{\Y}{ {\rm Y}}
\newcommand{\Z}{{\rm Z}}
\newcommand{\Ac}{\mathcal{A}}
\newcommand{\achi}{\Ac_{\chi}}
\newcommand{\bachi}{\overline\Ac_{\chi}}
\newcommand{\Lc}{\mathcal{L}}
\newcommand{\ol}{\overline}
\newcommand{\Pc}{\mathcal{P}}
\newcommand{\Qc}{\mathcal{Q}}
\newcommand{\Tc}{\mathcal{T}}
\newcommand{\Sc}{\mathcal{S}}
\newcommand{\Bc}{\mathcal{B}}
\newcommand{\Dc}{\mathcal{D}}
\newcommand{\Ec}{\mathcal{E}}
\newcommand{\Fc}{\mathcal{F}}
\newcommand{\Hc}{\mathcal{H}}
\newcommand{\Uc}{\mathcal{U}}
\newcommand{\Vc}{\mathcal{V}}
\newcommand{\Wc}{\mathcal{W}}
\newcommand{\Yc}{\mathcal{Y}}
\newcommand{\M}{\mathcal{M}}
\newcommand{\Ar}{{\rm A}}
\newcommand{\Br}{{\rm B}}
\newcommand{\Ir}{{\rm I}}
\newcommand{\Fr}{{\rm F}}
\newcommand{\Jr}{{\rm J}}
\newcommand{\Or}{{\rm O}}
\newcommand{\GL}{{\rm GL}}
\newcommand{\Spr}{{\rm Sp}}
\newcommand{\Rr}{{\rm R}}
\newcommand{\Zr}{{\rm Z}}
\newcommand{\gl}{\mathfrak{gl}}
\newcommand{\middd}{{\rm mid}}
\newcommand{\ev}{{\rm ev}}
\newcommand{\Pf}{{\rm Pf}}
\newcommand{\Norm}{{\rm Norm\tss}}
\newcommand{\oa}{\mathfrak{o}}
\newcommand{\spa}{\mathfrak{sp}}
\newcommand{\osp}{\mathfrak{osp}}
\newcommand{\g}{\mathfrak{g}}
\newcommand{\h}{\mathfrak h}
\newcommand{\n}{\mathfrak n}
\newcommand{\z}{\mathfrak{z}}
\newcommand{\Zgot}{\mathfrak{Z}}
\newcommand{\p}{\mathfrak{p}}
\newcommand{\sll}{\mathfrak{sl}}
\newcommand{\agot}{\mathfrak{a}}
\newcommand{\qdet}{ {\rm qdet}\ts}
\newcommand{\Ber}{ {\rm Ber}\ts}
\newcommand{\HC}{ {\mathcal HC}}
\newcommand{\cdet}{ {\rm cdet}}
\newcommand{\tr}{ {\rm tr}}
\newcommand{\gr}{ {\rm gr}}
\newcommand{\str}{ {\rm str}}
\newcommand{\loc}{{\rm loc}}
\newcommand{\Gr}{{\rm G}}
\newcommand{\sgn}{ {\rm sgn}\ts}
\newcommand{\ba}{\bar{a}}
\newcommand{\bb}{\bar{b}}
\newcommand{\bi}{\bar{\imath}}
\newcommand{\bj}{\bar{\jmath}}
\newcommand{\bk}{\bar{k}}
\newcommand{\bl}{\bar{l}}
\newcommand{\hb}{\mathbf{h}}
\newcommand{\Sym}{\mathfrak S}
\newcommand{\fand}{\quad\text{and}\quad}
\newcommand{\Fand}{\qquad\text{and}\qquad}
\newcommand{\For}{\qquad\text{or}\qquad}
\newcommand{\OR}{\qquad\text{or}\qquad}

\renewcommand{\theequation}{\arabic{section}.\arabic{equation}}

\newtheorem{thm}{Theorem}[section]
\newtheorem{lem}[thm]{Lemma}
\newtheorem{prop}[thm]{Proposition}
\newtheorem{cor}[thm]{Corollary}
\newtheorem{conj}[thm]{Conjecture}
\newtheorem*{mthm}{Main Theorem}
\newtheorem*{mthma}{Theorem A}
\newtheorem*{mthmb}{Theorem B}

\theoremstyle{definition}
\newtheorem{defin}[thm]{Definition}

\theoremstyle{remark}
\newtheorem{remark}[thm]{Remark}
\newtheorem{example}[thm]{Example}

\newcommand{\bth}{\begin{thm}}
\renewcommand{\eth}{\end{thm}}
\newcommand{\bpr}{\begin{prop}}
\newcommand{\epr}{\end{prop}}
\newcommand{\ble}{\begin{lem}}
\newcommand{\ele}{\end{lem}}
\newcommand{\bco}{\begin{cor}}
\newcommand{\eco}{\end{cor}}
\newcommand{\bde}{\begin{defin}}
\newcommand{\ede}{\end{defin}}
\newcommand{\bex}{\begin{example}}
\newcommand{\eex}{\end{example}}
\newcommand{\bre}{\begin{remark}}
\newcommand{\ere}{\end{remark}}
\newcommand{\bcj}{\begin{conj}}
\newcommand{\ecj}{\end{conj}}

\newcommand{\bal}{\begin{aligned}}
\newcommand{\eal}{\end{aligned}}
\newcommand{\beq}{\begin{equation}}
\newcommand{\eeq}{\end{equation}}
\newcommand{\ben}{\begin{equation*}}
\newcommand{\een}{\end{equation*}}

\newcommand{\bpf}{\begin{proof}}
\newcommand{\epf}{\end{proof}}
\newcommand{\whg}{\widehat{\mathfrak{g}}}

\def\beql#1{\begin{equation}\label{#1}}

\title{\Large\bf Quantum Sugawara operators in $\mathrm{U}_q(\widehat{\gl}_{M|N})$   }

\author{{Naihuan Jing,\ \ Ming Liu,\ \  Jian Zhang\ \  }}

\date{}  
\maketitle

\begin{abstract}
We construct Sugawara operators for the quantum affine superalgebra $\mathrm{U}_q(\widehat{\gl}_{M|N})$  in an explicit form which generalizes the results given in \cite{JLM} for the quantum affine algebras. We also calculate the Harish-Chandra images of
the Sugawara operators.

\end{abstract}

\vspace{5 mm}

\section{Introduction}

The Sugawara operators associated with a simple Lie algebra $\mathfrak{g}$ form a distinguished family of central elements in a suitable completion of the universal enveloping algebra of the corresponding affine Lie algebra at the critical level. These operators play a fundamental role in representation theory and mathematical physics, and admit deep connections with vertex algebras and the geometry of opers.
The remarkable properties of Sugawara operators have motivated extensive studies of their algebraic and geometric aspects, as well as  quantum analogues \cite{F, M}.

Quantum analogues of Sugawara operators arise in the framework of quantum affine algebras, where they appear as central elements of a completed algebra at the critical level. Their construction is intimately related to the R-matrix formalism and the associated $L$-operators satisfying the $RLL$ relations \cite{DE, RS}. In type A, this framework admits explicit constructions of quantum Sugawara operators, giving rise to a  family of central elements.
A systematic construction of these operators was initiated in the case of one-column Young diagrams \cite{fjmr:hs}, where explicit formulas for the generating series $S_{n}(z)$ were obtained, whose coefficients are quantum Sugawara operators. These elements belong to the center of an appropriate completion of the quantum affine algebra $\U_q(\widehat{\mathfrak{gl}}_n)$ at the critical level. This construction was subsequently generalized in~\cite{JLM} to arbitrary Young diagrams $\lambda$ containing at most $n$ rows.
The resulting central elements provide an explicit description of the center at the critical level and reveal remarkable combinatorial structures encoded by Young diagrams. The generalization can be viewed as a lift from algebraic generators towards a linear basis of the center.

The extension of these constructions to the superalgebra is a natural problem. However, \(\mathbb{Z}_2\)-grading significantly complicates the construction of central elements.
The usual arguments require substantial modifications, and new techniques are needed to treat the super case.
The super Yangian $\Y(\gl_{M|N})$
of the general linear Lie superalgebra $\gl_{M|N}$
was introduced by Nazarov \cite{Na, Na2}, who also proved that it contains a certain remarkable family of central elements. These elements arise as the coefficients of  quantum Berezinian, which is a super analogue of the quantum determinant for the Yangian $\Y(\gl_N)$ \cite{Mol, MNO}. It was conjectured by Nazarov \cite{Na} and later proved by Gow\cite{Gow} that the coefficients of Berezinian generate the center of $\Y(\gl_{M|N})$.
Bagnoli and Kozic \cite{BK} introduced the quantum Berezinian for double Yangian $\mathrm{DY}(\gl_{M|N})$ and proved that oefficients of
 quantum Berezinian
 are algebraically independent topological generators of the center of double Yangian $\mathrm{DY}(\gl_{M|N})$.
Molev and Ragoucy \cite{MR} used the Berezinian to construct
higher order Sugawara operators for the affine Lie superalgebra $\widehat{\mathfrak{gl}}_{M|N}$.

In parallel,  the quantum Berezinian for the quantum affine superalgebra $\U_q(\widehat{\mathfrak{gl}}_{M|N})$ was introduced in our previous work~\cite{JLZ}  and it was shown that its coefficients belong to the center of $\U_q(\widehat{\mathfrak{gl}}_{M|N})$. Another family of central elements was also constructed there in terms of the quantum Berezinian via a Liouville-type theorem.
These results provide super analogues of several classical constructions and provide important information about the center of quantum affine superalgebras. Nevertheless, In contrast to the central elements at arbitrary level, the critical level exhibits a much larger center.

The purpose of the present paper is to study the center of $\U_q(\widehat{\mathfrak{gl}}_{M|N})$ at the critical level and to give explicit constructions of general quantum Sugawara operators for quantum affine superalgebras $\U_q(\widehat{\mathfrak{gl}}_{M|N})$. Working within the $RLL$ formalism, we construct explicit generating series whose coefficients belong to the center of a suitable completion of $\U_q(\widehat{\mathfrak{gl}}_{M|N})$ at the critical level. These series can be regarded as super analogues of the quantum Sugawara operators in $\U_q(\widehat{\mathfrak{gl}}_n)$. Our construction builds on the structure of the $L$-operators and employs the approach of~\cite{JLM}, while also reflecting new features specific to the superalgebra setting. We further investigate the Harish-Chandra image of these operators.
 The Harish-Chandra image provides an effective description of the resulting Sugawara operators in terms of commutative variables and offers valuable insight into the structure of the center at the critical level. We expect that these results will contribute to a deeper understanding of the representation theory of quantum affine superalgebras and related problems in quantum integrable systems.

\section{Quantum affine superalgebra $\mathrm{U}_q(\widehat{\gl}_{M|N})$}

In this section, we present the RLL realization of the quantum affine superalgebra $\mathrm{U}_q(\widehat{\gl}_{M|N})$.
For any $1\leq i\leq M+N$, the parity of $i$ is defined as
\[
\bar{i}=\left\{\begin{array}{cc}
0&\text{ if }i\leq M\\
1&\text{ if }i> M,
\end{array}\right.
\]
Let $q_i=q^{1-2\bar i}$ for $i\in \{1,\ldots,M+N\}$ ,$E_{ij}$ be the unit matrix in $\End \mathbb{C}^{M|N} $,
the $R$-matrix $R(z,w) \in (\End  \mathbb{C}^{M|N} \otimes \End  \mathbb{C}^{M|N} )[z,w]$ is defined as
\begin{eqnarray}
\begin{array}{rcl}
R(z,w) &=&  \sum\limits_{i\in I}(zq_i - wq_i^{-1}) E_{ii} \otimes E_{ii}  + (z-w) \sum\limits_{i \neq j} E_{ii} \otimes E_{jj} \\
&\ & + z\sum\limits_{i<j}(q_j-q_j^{-1})  E_{ij} \otimes E_{ji}+w \sum\limits_{i>j} (q_j-q_j^{-1}) E_{ij} \otimes E_{ji} .
\end{array}
\end{eqnarray}
It satisfies the following {\it quantum Yang-Baxter equation}
\begin{equation}\label{QYBE}
R_{12}(z_1,z_2)R_{13}(z_1,z_3)R_{23}(z_2,z_3) = R_{23}(z_2,z_3)R_{13}(z_1,z_3)R_{12}(z_1,z_2).
\end{equation}
Introduce the graded permutation operator $P$ on the tensor product $\mathbb{C}^{M|N}  \otimes \mathbb{C}^{M|N}  $ such that
$$P=\sum_{i,j=1}^{M+N}(-1)^{\bar{j}}\;E_{ij}\otimes E_{ji}.$$ Let $R = R(1,0),R' = P R^{-1}P$. Then $R(z,w) = z R - w R'$ and $R-R' =(q - q^{-1}) P$.
We also need the  $R$ matrix:
\beq
\ol R(z/w)=\frac{R(z,w)}{zq-wq^{-1}}.
\eeq
Denote $\ol R_{21}(z)=P_{12}\ol R_{12}(z)P_{12}$, then
\beql{Rin}
\ol R_{12}({z\over w})\ol R_{21}({w\over z})=1,\quad
\ol R_{21}({w\over z}) = \ol R_{q^{-1}}({z\over w}).
\eeq

For  a matrix $A=\sum\limits _{i,j=1}^{M+N}a_{ij}E_{ij}$,   the supertranspose is defined by
$$A^{st}=\sum_{i,j}(-1)^{\bar{i}(\bar{i}+\bar{j})}a_{ji}E_{ij},$$
and the supertrace $str$ by
\[
str(A)=\sum_{i=1}^{M+N}(-1)^{\bar{i}}a_{ii}.
\]
For any  $a\in\{1,2,\ldots,k\}$ we will denote by $st_a$ the corresponding partial transposition on the
algebra $(\End \mathbb C^{M|N})^{\otimes k}$ which acts as $st$ on the $a$-th copy of $\End \mathbb C^{M|N}$ and as the identity map on all
the other tensor factors.

For $M\neq N$, we introduce the following normalized $R$-matrix
\beql{rf}
R(x)=f(x)\tss \overline R(x),
\eeq
where
\ben
f(x)=1+\sum_{k=1}^{\infty}f_kx^k,\qquad f_k=f_k(q),
\een
is a formal power series in $x$ 
uniquely determined by the relation
\ben
f(xq^{2N-2M})=f(x)\ts\frac{(1-xq^{-2})\tss(1-xq^{2N-2M+2})}{(1-x)\tss(1-xq^{2N-2M})}.
\een

The $R$-matrix \eqref{rf} satisfies the {\em crossing symmetry relations} \cite{GZ}:
\beql{cs}
\big(R_{12}(x)^{-1}\big)^{st_2} D_2 R_{12}(xq^{2N-2M})^{st_2}=D_2
 \text{ and }
R_{12}(xq^{2N-2M})^{st_1}\tss D_1\big(R_{12}(x)^{-1}\big)^{st_1}=D_1,
\eeq
where $D$ denotes the diagonal $n\times n$ matrix
\beql{d}
D=\diag\big[q^{2},q^{4},\ldots,q^{2M},q^{2M},\ldots,q^{2M-2N+2}\big].
\eeq

Following
\cite{GZ,Z},
We introduce the quantum affine superalgebra in the R-matrix format.
\defin
The quantum affine superalgebra $\mathrm U_q(\wh\gl_{M|N})$ ,$M\neq N$, is the superalgebra
generated by ${l^{\pm}_{ij}}^{(r)}$  where $1\leq i, j\leq M+N$ and $r$ runs over nonnegative integers. Let $L^{\pm}(u)=(l_{ij}^{\pm}(u))$ be the matrix
\begin{equation}
\begin{split}
L^{\pm}(u)= \sum_{i,j=1}^{M+N}l^{\pm}_{ij}(u)\otimes E_{ij},
\end{split}
\end{equation}
where $l^{\pm}_{ij}(u)$  are formal series in $u^{\pm 1}$ respectively:
\begin{equation}\label{l series}
\begin{split}
l^{\pm}_{ij}(u)= \sum_{r=0}^{\infty} {l^{\pm}_{ij}}^{(r)}u^{\pm r}.
\end{split}
\end{equation}
The defining relations are
\begin{align} \label{RLL-0}
{l^-_{ji}}^{(0)}={l^+_{ij}}^{(0)}=0,\ 1\leq i<j\leq M+N, \\ \label{RLL-1}
l{^-_{ii}}^{(0)} {l^+_{ii}}^{(0)}={l^+_{ii}}^{(0)}{l^-_{ii}}^{(0)}=1,\ 1\leq i\leq M+N, \\ \label{RLL}
R(z/w)L^{\pm}_{1} (z) L^{\pm}_{2} (w)= L^{\pm}_{2} (w) L^{\pm}_{1} (z)R(z/w),\\ \label{RLL cros}
R(z q^{c} /w) L^{+}_{1} (z) L^{-}_{2} (w)= L^{-}_{2} (w) L^{+}_{1} (z)R(z q^{-c}  /w).
\end{align}

Note that when $N=0$ the quantum affine superalgebra $\mathrm U_q(\wh\gl_{M|0})$ descends to the quantum affine algebra
$\mathrm U_q(\wh\gl_{M})$ studied by Ding and Frenkel \cite{DF}.

Denote by $\U_q(\wh\gl_{M|N})_{\text{\rm cri}}$ the quantum affine algebra
{\em at the critical level $c=N-M$}, which is
the quotient of $\U_q(\wh\gl_{M|N})$ by the relation $q^c=q^{N-M}$.
Its completion $\wt\U_q(\wh\gl_{M|N})_{\text{\rm cri}}$ is defined as
the inverse limit
\beql{compl}
\wt\U_q(\wh\gl_{M|N})_{\text{\rm cri}}=\lim_{\longleftarrow}
\U_q(\wh\gl_{M|N})_{\text{\rm cri}}/J_p, \qquad p>0,
\eeq
where $J_p$ denotes the left ideal of $\U_q(\wh\gl_{M|N})_{\text{\rm cri}}$ generated by all elements
$l^{-}_{ij}[r]$ with $r\geqslant p$.
Elements of the center $\Zr_q(\wh\gl_{M|N})$ of $\wt\U_q(\wh\gl_{M|N})_{\text{\rm cri}}$
are known as ({\em quantum}) {\em Sugawara operators}.

\section{Hecke algebra}
The {\em Hecke algebra} $\mathcal{H}_m$ of type $A_{m-1}$ is generated by $T_1,\dots,T_{m-1}$ subject to the well‑known relations:
	\[
	(T_i - q)(T_i + q^{-1}) = 0,\qquad
	T_i T_{i+1} T_i = T_{i+1} T_i T_{i+1},\qquad
	T_i T_j = T_j T_i\;\;(|i-j|>1).
	\]
	For any permutation $\sigma \in \mathfrak{S}_m$, pick a reduced decomposition $\sigma = \sigma_{i_1}\cdots\sigma_{i_l}$ where each $\sigma_i = (i,i+1)$ is an adjacent transposition, and set $T_\sigma = T_{i_1}\cdots T_{i_l}$. The defining relations of $\mathcal{H}_m$ guarantee that $T_\sigma$ does not depend on the chosen reduced decomposition. The collection $\{T_\sigma\}_{\sigma\in\mathfrak{S}_m}$ forms a linear basis of $\mathcal{H}_m$.
	An involutive anti‑automorphism $*$ is defined on $\mathcal{H}_m$ by $T_\sigma^* = T_{\sigma^{-1}}$.

\paragraph{Young basis.}

	Let $\lambda = (\lambda_1,\dots,\lambda_\ell)$ be a partition of $m$, identified with its Young diagram. Its conjugate partition is denoted $\lambda' = (\lambda_1',\dots,\lambda_r')$, where $\lambda_j'$ counts the number of boxes in column $j$. For a box $\alpha = (i,j)\in \lambda$, the hook length is $h(\alpha) = \lambda_i + \lambda_j' - i - j + 1$, and the content is $c(\alpha) = j - i$.
	
	A $\lambda$-tableau $\mathcal{U}$ is a filling of the diagram $\lambda$ with numbers taken from $\{1,\dots,n\}$. It is called \emph{standard} if the entries increase strictly along each row and down each column; it is \emph{semistandard} if the entries increase weakly along rows and strictly down columns.
	Given a standard $\lambda$-tableau $\Lambda$, write $c_k(\Lambda)$ for the content of the cell containing the entry $k$, and define $d_k(\Lambda) = c_{k+1}(\Lambda) - c_k(\Lambda)$.
	
	The irreducible representations of $\mathcal{H}_m$ over $\mathbb{C}$ are labelled by partitions of $m$. For a partition $\lambda \vdash m$, let $V_\lambda$ be the corresponding irreducible module and let
	\beql{phila}
    \vp_{\la}:\Hc_{m}\to \End V_{\la}
    \eeq
	be the associated algebra homomorphism. The space $V_\lambda$ carries an $\mathcal{H}_m$-invariant inner product $\langle \cdot,\cdot\rangle = \langle \cdot,\cdot\rangle_{V_\lambda}$ satisfying $\langle h v_1, v_2\rangle = \langle v_1, h^* v_2\rangle$ for all $h\in\mathcal{H}_m$ and $v_1,v_2\in V_\lambda$ ( see \cite[Sec.~3.1]{dj:bi}.).
	
	There exists an orthonormal \emph{Young basis} $\{v_\Lambda\}$ of $V_\lambda$, indexed by the standard $\lambda$-tableaux $\Lambda$. Denote $f_\lambda = \dim V_\lambda$, which equals the number of such standard tableaux. Using the seminormal form developed in \cite{dj:bi} and \cite{h:rh}, the action of the generators $T_k$ on this basis can be made completely explicit. After a suitable normalization, for any $k \in \{1,\dots,m-1\}$ one obtains
	\beql{eq:Young basis}
    T_{k} \ts v^{}_{\La}=\frac{q^{d_k(\La)}}{[d_k(\La)]_q} v^{}_{\La}+\sqrt{1-\frac{1}{[d_k(\La)]_q^2}}
    \ts v^{}_{\si_{k} \La},
    \eeq
	where $[n]_q = \frac{q^n - q^{-n}}{q - q^{-1}}$, the tableau $\sigma_k\Lambda$ results from $\Lambda$ by swapping the entries $k$ and $k+1$, and $v_{\sigma_k\Lambda}=0$ whenever $\sigma_k\Lambda$ fails to be standard.
	
	For any skew diagram $\theta$ consisting of $m$ boxes, the same formulas \eqref{eq:Young basis} define a representation of $\mathcal{H}_m$ on the linear span $V_\theta$ of the vectors $v_\Lambda$ labelled by standard tableaux $\Lambda$ of shape $\theta$.
	
\paragraph{Symmetrizing trace and Schur elements.}
Define the {\em symmetrizing trace} on the Hecke algebra as the following linear map
\ben
\tau:\Hc_{m}\to\CC,\qquad \tau(T_{\si})=\de_{\si,e},
\een
where $e$ denotes the identity element of $\Sym_m$. This induces a {\em symmetric algebra}
structure on $\Hc_{m}$ associated with the following bilinear form
\ben
\Hc_{m}\ot\Hc_m\to\CC,\qquad h_1\ot h_2\mapsto\tau(h_1h_2),
\een
which is symmetric and non-degenerate.  The dual basis of
$\{T_{\si}\}_{\si \in \Sym_m}$ with respect to the form is
$\{ T_{\si^{-1}}\}_{\si \in \Sym_m}$ \cite{GP, m:rt}.

For any $\pi\in\End V_{\la}$, define $I(\pi)\in \End V_{\la}$ by
\ben
I(\pi)=\sum_{\si\in \Sym_m}\vp_\la(T_{\si})\ts \pi\ts\vp_\la(T_{\si^{-1}}).
\een
	A direct verification shows that $I(\pi)$ actually lies in $\End_{\mathcal{H}_m} V_\lambda$; therefore Schur's lemma forces it to be a scalar multiple of the identity. Consequently we obtain the key relation
	\beql{iuschur}
I(\pi)=c_{\la}\tss\tr^{}_{V_{\la}}(\pi)\tss{\rm id}^{}_{V_{\la}},
\eeq
	which will be essential in the proof of  Lemma \ref{lem:ASchur}. Here $c_\lambda$ is the \emph{Schur element} given by the Steinberg formula \cite[Thm~4.64]{m:rt}:
	\ben
c_{\la}=\prod_{\al\in\la} q^{c(\al)}\tss [h(\al)]_q.
\een
	Let $\chi_\lambda$ denote the character of the $\mathcal{H}_m$-module $V_\lambda$. Then the symmetrizing trace decomposes as
	\beql{tauchi}
\tau=\sum_{\la\vdash m}\frac{\chi_{\la}}{c_{\la}}.
\eeq

\paragraph{Idempotents and matrix units.}
The Hecke algebra $\mathcal{H}_m$ is semisimple and admits an isomorphism
\beql{isomhecke}
\Hc_m \cong \bigoplus_{\la \vdash m} \Mat_{f_{\la}}(\CC).
\eeq
Via this isomorphism, the matrix units $e^{\la}_{\La,\La'} \in \Mat_{f_{\la}}(\CC)$ corresponding to $V_\lambda$
can be identified with elements of $\Hc_m$ by
\beql{mat-units}
e^{\la}_{\La,\La'}=\frac{1}{c_{\la}}\sum_{\si\in \Sym_m}
\langle T_{\si^{-1}}v^{}_{\La},v^{}_{\La'}\rangle^{}_{V_{\la}}T_{\si}.
\eeq
	The characters satisfy the orthogonality relations: for any partitions $\lambda,\mu\vdash m$,
\beql{orthcha}
\sum_{\si\in \Sym_m}\chi_{\la}(T_{\si})\chi_{\mu}(T_{\si^{-1}})=\de_{\la,\mu}c_{\la}f_{\la}.
\eeq

The diagonal matrix units $e^{\la}_{\La}=e^{\la}_{\La \La}$ are primitive idempotents of $\Hc_m$.
Note the immediate properties of the primitive idempotents:
\ben
e^{\la}_{\La} e^{\la}_{\Ga}=0 \quad \text { if } \quad \La \neq \Ga\Fand
(e^{\la}_{\La})^{2}=e^{\la}_{\La}
\een
together with the decomposition of the identity in $\Hc_m$
\ben
1=\sum_{\la \vdash m} \sum_{\sh(\La)=\la} e^{\la}_{\La}.
\een

The following lemma given in \cite{JLM} will play an important role in the proof of Lemma \ref{lem:RT}.
\ble\label{lem:et}
Let $\La$ be a standard $\la$-tableau and let  $k \in\{1, \dots, m-1\}$.
Then
\ben
e^{\la}_{\La}\left(T_{k}-\frac{q^{d_k(\La)}}{[d_k(\La)]_q}\right)=\sqrt{1-\frac{1}{[d_k(\La)]_q^2}}
\ts e^{\la}_{\La,\si_{k} \La}
\een
and
\ben
\left(T_{k}+\frac{q^{-d_k(\La)}}{[d_k(\La)]_q}\right)e^{\la}_{\si_{k} \La}
=\sqrt{1-\frac{1}{[d_k(\La)]_q^2}} \ts e^{\la}_{\La,\si_{k} \La},
\een
assuming that $e^{\la}_{\La,\si_{k} \La}=0$ if the tableau $\si_{k} \La$ is not standard.
\ele
\section{Quantum Sugawara operators}
In this section, we use the Fusion procedure \cite{MI, Naz} to construct quantum Sugawara operators for $\mathrm{U}_q(\widehat{\gl}_{M|N})$.

The $R$-matrix $\ol  R(z)$ can be written as
\ben
\ol R(z)=\frac{z-1}{qz-q^{-1} }\ts\Big(R+\frac{q-q^{-1}}{z-1}P\Big),
\een
where
\begin{eqnarray}
\begin{array}{rcl}
R &=&  \sum\limits_{i\in I} q_i   E_{ii} \otimes E_{ii}  +   \sum\limits_{i \neq j} E_{ii} \otimes E_{jj}
+ \sum\limits_{i<j}(q_j-q_j^{-1})  E_{ij} \otimes E_{ji}.
\end{array}
\end{eqnarray}
Setting $\check{R}=PR$, we get
\ben
\check{R}_{k}\check{R}_{k+1}\check{R}_{k}=\check{R}_{k+1}\check{R}_{k}\check{R}_{k+1}\Fand
(\check{R}_{k}-q)(\check{R}_{k}+q^{-1})=0,
\een
where $\check{R}_{k}=P_{k,k+1}R_{k,k+1}$.

For each $k=1,\dots, m-1$, introduce the $\Hc_{m}$-valued rational functions in two variables
$x, y$ by
\ben
T_{k}(x,y)=T_k+\frac{q-q^{-1}}{xy^{-1}-1},
\een
which satisfy the relations
\ben
T_{k}(x,y)\ts T_{k+1}(x,z)\ts T_{k}(y,z)=T_{k+1}(y,z)\ts T_{k}(x,z)\ts T_{k+1}(x,y)
\een
and
\ben
T_{k}(x, y)\ts T_{k}(y, x)=1-\frac{\left(q-q^{-1}\right)^{2} x y}{(x-y)^{2}}.
\een

There exsits a representation of the Hecke algebra $\Hc_{m}$ on the tensor product space
$(\CC^{M|N})^{\ot m}$
defined by
\beql{haact}
T_{k}\mapsto \check{R}_{k},\qquad k=1,\dots,m-1.
\eeq
Furthermore, we find that under this action of the Hecke algebra,
\ben
T_{k}(x,y)\mapsto \check{R}_{k}(x/y)\qquad\text{with}\quad
\check{R}_{k}(z)=\check{R}_k+\frac{q-q^{-1}}{z-1}.
\een

Let $\La$ be a standard tableau of shape $\la$. We will keep the notation $c_{k}(\La)$
for the content $j-i$ of the box $(i,j)$ of $\la$ occupied by $k$ in $\La$.
Equip the set of all pairs $(i, j)$
with $1 \leqslant i<j \leqslant m$ with the following ordering. The pair $(i, j)$ precedes
$(i^{\tss\prime}, j^{\tss\prime})$ if $j<j^{\tss\prime}$, or if $j=j^{\tss\prime}$
but $i<i^{\tss\prime}$. Set
\beql{F}
T_{\La}(z_{1}, \dots, z_{m})=
\overrightarrow{\prod\limits_{(i, j)}}\ts
T_{j-i}(z_{i},   z_{j})
\eeq
with the ordered product taken over the set of pairs. This is a rational function
in variables $z_{1}, \dots, z_{m}$ taking values
in $\Hc_{m}$.

For any $k=1,\ldots,m$, let $w_k$ denote the unique longest element of the symmetric
group $S_k$ which is regarded as the natural subgroup of $S_m$. The corresponding
elements  $T_{w_k}\in\Hc_{m}$
are   given by $t_{w_1}=1$ and
\beq
\bal
T_{w_k}=&T_1(T_2T_1)\cdots (T_{k-2}\cdots T_1) (T_{k-1}T_{k-2}\cdots T_1)\\
=&(T_1\cdots T_{k-1})(T_1\cdots T_{k-2})\cdots  (T_1 T_2)T_1, \qquad k=2,\ldots,m .
\eal
\eeq
According to \cite[Thm~3.3]{imo:ih},
the primitive idempotents $e^{\la}_{\La}$ can be obtained by
the consecutive evaluations
\beql{theor}
e^{\la}_{\La} = \frac{1}{c_{\la'}}\ts T_{\La}(z_{1}, \dots, z_{m})\ts T_{w_m}^{-1}\ts
\big|_{z_1=q^{-2 c_{1}(\La)} }\big|_{z_2=q^{-2 c_{2}(\La)} } \dots \big|_{z_m=q^{-2 c_{m}(\La)} },
\eeq
which are all well-defined; the rational functions are regular at
the evaluation points at each step, and $c_{\la'}$ is the Schur element associated with
 the conjugate partition $\la^{\prime}$.

Given a standard tableau $\La$ of shape $\la\vdash m$,
we denote by $\check{R}_{\La}(z_{1}, \dots, z_{m})$,  $\check{R}_{w_m}$ , $\Ec^{\la}_{\La}$ (or just $\Ec_{\La}$) the image of the image of
$T_{\La}(z_{1}, \dots, z_{m})$ , $T_{w_m}$ and $e^{\la}_{\La}$
under the action \eqref{haact}
respectively. Clearly, $\Ec_{\La}^2=\Ec_{\La}$, and the fusion formula
\eqref{theor} implies the relation
\beql{fusr}
\Ec_{\La} = \frac{1}{c_{\la'}}\ts \check{R}_{\La}(z_{1}, \dots, z_{m})\ts \check{R}_{w_m}^{-1}\ts
\big|_{z_1=q^{-2 c_{1}(\La)} }\big|_{z_2=q^{-2 c_{2}(\La)} } \dots \big|_{z_m=q^{-2 c_{m}(\La)} }.
\eeq

Set
\ben
L^{\pm}_{\La}(z)=L^\pm_1(zq^{2c_1(\La)})\dots L^\pm_m(zq^{2c_{m}(\La)})
\een
and $L_{\La}(z)=L^{+}_{\La}(z)L^{-}_{\La}(zq^{M-N})^{-1}$.

\ble\label{lem:LE}
We have the relations
\ben
L^{+}_{\La}(z)\tss\Ec_{\La}=\Ec_{\La}\tss L^{+}_{\La}(z)\tss\Ec_{\La},
\een
\ben
L^{-}_{\La}(z )^{-1}\Ec_{\La}=\Ec_{\La}\tss L^{-}_{\La}(z )^{-1}\Ec_{\La},
\een
and hence
\ben
L_{\La}(z)\tss \Ec_{\La}=\Ec_{\La}\tss L_{\La}(z)\tss\Ec_{\La}.
\een
\ele

\bpf
It follows from the RTT relations that
\beql{eqRLL}
\bal
\check{R}_{\La}(z_{1}, \dots, z_{m})
&L^{\pm}_1(z/z_m)\dots L^{\pm}_m(z/z_1) \\
=&L^{\pm}_1(z/z_1)\dots L^{\pm}_m(z/z_m) \check{R}_{\La}(z_{1}, \dots, z_{m}).
\eal
\eeq
Multiply both sides by the element $\check{R}_{w_m}^{-1}$ from the right and apply consecutive
evaluations by setting $z_1=q^{-2 c_{1}(\La)} ,z_2=q^{-2 c_{2}(\La)} , \dots ,z_m=q^{-2 c_{m}(\La)}$ to derive from \eqref{fusr} that
\ben
\Ec_{\La}\check{R}_{w_m}\tss L_{1}^{+}(zq^{2c_m(\La)})\dots
L_{m}^{+}(zq^{2c_1(\La)})\check{R}_{w_m}^{-1}
=L_{1}^{+}(zq^{2c_1(\La)})\dots L_{m}^{+}(zq^{2c_m(\La)})\ts\Ec_{\La}.
\een
This shows that the right hand side stays unchanged when it is multiplied by $\Ec_{\La}$
from the left thus proving the first relation.

To prove the second relation, we   multiply both sides of \eqref{eqRLL} by the inverses of the products of the $L$-matrices
to get
\ben
\bal
\check{R}_{\La}(z_{1}, \dots, z_{m})\ts L_{m}^{-}&(z/z_1)^{-1}\dots
L_{1}^{-}(z/z_m )^{-1}\\
&=
L_{m}^{-}(z/z_m)^{-1}\dots L_{1}^{-}(z/z_1)^{-1}
\ts\check{R}_{\La}(z_{1}, \dots, z_{m}).
\eal
\een
Then the second relation can be proved
by the same argument as for the first relation.
\epf
The Yang-Baxter equation implies that
\beq
\check{R}_{ij}(z_i/z_j)R_{0j}(zz_i)R_{0i}(zz_j)=R_{0j}(zz_j)R_{0i}(zz_i)\check{R}_{ij}(z_i/z_j).
\eeq

Therefore, by replacing the $L$-matrices with the elements $R_{0i}(z/z_i)$ in the proof
of Lemma~\ref{lem:LE}, we get its following counterpart.

\ble\label{lem:RE}
We have the relations
\ben
\Ec_{\La}R_{0m}\Big( zq^{-2c_m(\La)} \Big)\dots R_{01}
\Big({zq^{-2c_1(\La)}}\Big)\Ec_{\La}=
R_{0m}\Big( zq^{-2c_m(\La)} \Big)\dots R_{01}
\Big({zq^{-2c_1(\La)}}\Big)\Ec_{\La}
\een
and
\ben
\Ec_{\La}R_{01}
\Big({zq^{-2c_1(\La)}}\Big)^{-1}\dots
R_{0m}\Big( zq^{-2c_m(\La)} \Big)^{-1}\Ec_{\La}=
R_{01}
\Big({zq^{-2c_1(\La)}}\Big)^{-1}\dots
R_{0m}\Big( zq^{-2c_m(\La)} \Big)^{-1}\Ec_{\La}.
\een
\ele

For each $a=1,\dots,k$ the supertrace $\str_a$
with respect to the $a$-th copy of $\End\CC^{M|N}$
in $\big(\End\CC^{M|N}\big)^{\ot\ts k}$ is the linear map
\ben
\str_a: U_q(\whg)\ot \big(\End\CC^{m|n}\big)^{\ot\ts k}
\to U_q(\whg)\ot\big(\End\CC^{m|n}\big)^{\ot\ts k-1}
\een
defined by
\ben
\bal
\str_a:u\otimes x_1\ot\dots \ot x_{a-1}&\ot e_{ij}\ot
x_{a+1}\ot \dots\ot x_k \\
{}&{}\mapsto \de_{ij}\ts (-1)^{\bi} u\ot x_1\ot\dots \ot x_{a-1}\ot x_{a+1}
\ot \dots\ot x_k\ot y\ts.
\eal
\een

Introduce the Laurent series $S_{\La}(z)$ in $z$ by
\ben
S_{\La}(z)=\str^{}_{1,\dots,m}\ts L_{\La}(z)D_1\dots D_m \Ec_{\La},
\een
the trace is taken over all $m$ copies of $\End\CC^{M|N}$.
All coefficients of the series $S_{\La}(z)$
are elements of the algebra $\wt\U_q(\wh\gl_{M|N})_{\text{\rm cri}}$.

The following is our main result which provides explicit formulas
for quantum Sugawara operators.

\bth\label{thm:centr}
The coefficients of the series $S_{\La}(z)$ belong to the center $\Zr_q(\wh\gl_{M|N})$ of the completed quantum
affine superalgebra at the critical level $\wt\U_q(\wh\gl_{M|N})_{\text{\rm cri}}$. Moreover,
$S_{\La}(z)$ does not depend on the standard $\la$-tableau $\La$
and only depends on the Young diagram $\la$.
\eth

\bpf
It will be sufficient to verify what $S_{\La}(z)$ commutes with $L^{\pm}(w)$. The calculations are similar in both cases
and so we will only consider $L^+(z)$ and follow the arguments of \cite[Sec.~3]{fjmr:hs}.

The coefficient of $E_{ij}\otimes E_{kl}$ in the $R$-matrix $R(z,w)$ is nonzero only if
$\bar i+\bar j +\bar k+\bar l=0\ mod\ 2$. Then $L_0^+(w)\Ec_{\La}=\Ec_{\La}L_0^+(w)$,
It is obvious that
$L_0^+(w)D_i=D_iL_0^+(w)$.

By using \eqref{RLL} and \eqref{RLL cros}
we get
\ben
\bal
L^+_0(w)L_{\La}(z)D_1\dots &  D_m \Ec_{\La}
=
R_{01}\Big(\frac{w}{z q^{2c_1(\La)}} \Big)^{-1}\dots
R_{0m}\Big(\frac{w}{z q^{2c_m(\La)}}\Big)^{-1}L_{\La}(z)\\[0.2em]
&\times R_{0m}\Big(\frac{w}{zq^{2M-2N+2c_m(\La)}} \Big)\dots
R_{01}\Big(\frac{w}{zq^{2M-2N+2c_1(\La)}} \Big)
D_1\dots D_m\Ec_{\La}L^+_0(w).
\eal
\een
Therefore, to conclude that $L^+_0(z)\tss S_{\La}(w)=S_{\La}(w)\tss L^+_0(z)$ we need
to show that the supertrace
\beq\label{str}
\bal
\str^{}_{1,\dots,m}
&R_{01}\Big(\frac{w}{z q^{2c_1(\La)}} \Big)^{-1}\dots
R_{0m}\Big(\frac{w}{z q^{2c_m(\La)}}\Big)^{-1}L_{\La}(z)\\[0.2em]
&\times R_{0m}\Big(\frac{w}{zq^{2M-2N+2c_m(\La)}} \Big)\dots
R_{01}\Big(\frac{w}{zq^{2M-2N+2c_1(\La)}} \Big)
D_1\dots D_m\Ec_{\La}
\eal
\eeq
equals $S_{\La}(w)$.

Set
\ben
X= R_{01}\Big(\frac{w}{z q^{2c_1(\La)}} \Big)^{-1}\dots
R_{0m}\Big(\frac{w}{z q^{2c_m(\La)}}\Big)^{-1}L_{\La}(z)\Ec_{\La}
\een
and
\ben
Y=  R_{0m}\Big(\frac{w}{zq^{2M-2N+2c_m(\La)}} \Big)\dots
R_{01}\Big(\frac{w}{zq^{2M-2N+2c_1(\La)}} \Big)
D_1\dots D_m\Ec_{\La}.
\een
Then
\ben
\str^{}_{1,\dots,m}\ts XY=\str^{}_{1,\dots,m}\ts X^{st_1\dots st_m}Y^{st_1\dots st_m}.
\een
Using the relation $\Ec_{\La}^2=\Ec_{\La}$ and applying Lemma \ref{lem:LE} and \ref{lem:RE},  we have that
\ben
\bal
X= R_{01}\Big(\frac{w}{z q^{2c_1(\La)}} \Big)^{-1}\dots
R_{0m}\Big(\frac{w}{z q^{2c_m(\La)}}\Big)^{-1}L_{\La}(z)\Ec_{\La}\\
=\Ec_{\La} R_{01}\Big(\frac{w}{z q^{2c_1(\La)}} \Big)^{-1}\dots
R_{0m}\Big(\frac{w}{z q^{2c_m(\La)}}\Big)^{-1}  L_{\La}(z)\Ec_{\La}.
\eal
\een
Since $\Ec_{\La}$ and $R$  are even operators,
\ben
\bal
X^{st_1\dots st_m}= \Ec_{\La}^{st_1\dots st_m}L_{\La}^{st_1\dots st_m}(z)
\Bigg(R_{01}\Big(\frac{w} {z q^{2c_1(\La)}}\Big)^{-1}\Bigg)^{st_1}\dots
\Bigg(R_{0m}\Big(\frac {w} {z q^{2c_m(\La)}} \Big)^{-1}\Bigg)^{st_m}\\
=\Ec_{\La}^{st_1\dots st_m}L_{\La}^{st_1\dots st_m}(z)
\Bigg(R_{01}\Big(\frac{w} {z q^{2c_1(\La)}}\Big)^{-1}\Bigg)^{st_1}\dots
\Bigg(R_{0m}\Big(\frac {w} {z q^{2c_m(\La)}} \Big)^{-1}\Bigg)^{t_m} \Ec_{\La}^{st_1\dots st_m}.\\
Y^{st_1\dots st_m}=
\Ec_{\La}^{st_1\dots st_m} D_1\dots D_m\tss R_{0m}\Big(\frac{w}{zq^{2M-2N+2c_m(\La)}} \Big)^{st_m}\dots
R_{01}\Big(\frac{w}{zq^{2M-2N+2c_1(\La)}} \Big)^{st_1}
\eal
\een
Therefore, we get
\begin{multline}
\str^{}_{1,\dots,m}\ts X^{st_1\dots st_m}Y^{st_1\dots st_m}\\
{}=\str^{}_{1,\dots,m}\ts \Ec_{\La}^{st_1\dots st_m}L_{\La}^{st_1\dots st_m}(z)
\Bigg(R_{01}\Big(\frac{w} {z q^{2c_1(\La)}}\Big)^{-1}\Bigg)^{st_1}\dots
\Bigg(R_{0m}\Big(\frac {w} {z q^{2c_m(\La)}} \Big)^{-1}\Bigg)^{st_m}\\[0.5em]
\times  \Ec_{\La}^{st_1\dots st_m} D_1\dots D_m\tss R_{0m}\Big(\frac{w}{zq^{2M-2N+2c_m(\La)}} \Big)^{st_m}\dots
R_{01}\Big(\frac{w}{zq^{2M-2N+2c_1(\La)}} \Big)^{st_1}\\
=\str^{}_{1,\dots,m}\ts \Ec_{\La}^{st_1\dots st_m}L_{\La}^{st_1\dots st_m}(z)
\Bigg(R_{01}\Big(\frac{w} {z q^{2c_1(\La)}}\Big)^{-1}\Bigg)^{st_1}\dots
\Bigg(R_{0m}\Big(\frac {w} {z q^{2c_m(\La)}} \Big)^{-1}\Bigg)^{st_m}\\[0.5em]
\times    D_1\dots D_m\tss R_{0m}\Big(\frac{w}{zq^{2M-2N+2c_m(\La)}} \Big)^{st_m}\dots
R_{01}\Big(\frac{w}{zq^{2M-2N+2c_1(\La)}} \Big)^{st_1}.
\non
\end{multline}
Use the first crossing symmetry relation in \eqref{cs} to conclude that
\ben
\bal
\str^{}_{1,\dots,m}\ts X^{st_1\dots st_m}Y^{st_1\dots st_m}
{}&=\str^{}_{1,\dots,m}\ts\Ec_{\La}^{st_1\dots st_m}L_{\La}^{st_1\dots st_m}(z) D_1\dots D_m
\\
{}&=\str^{}_{1,\dots,m}\ts D_1\dots D_m\Ec_{\La}^{st_1\dots st_m}L_{\La}^{st_1\dots st_m}(w)
\\
{}&=\str^{}_{1,\dots,m}\ts L_{\La}(w)\Ec_{\La}
D_1\dots D_m
\eal
\een
which coincides with $S_{\La}(w)$. This proves the first part of the theorem.

To prove the second part of the theorem , it's sufficient to verify the relation
\beql{lalapr}
\str^{}_{1,\dots,m}\ts L_{\La}(z)D_1\dots D_m \Ec_{\La}
=\str^{}_{1,\dots,m}\ts L_{\La'}(z)D_1\dots D_m \Ec_{\La'},
\eeq
where $\La'=\si_k\La$ is also a standard tableau for $k\in\{1,\dots,m-1\}$.
Set $d_k=c_{k+1}(\La)-c_{k}(\La)$ and note that both operators
$\check{R}_k(q^{-2d_{k}})$ and $\check{R}_k(q^{2d_k})$ are invertible.
By taking into account the property
$P\overline R(x) P=\overline R(x^{-1})^{-1}$
of the $R$-matrix \eqref{Rin} and using \eqref{RLL} we get the relation
\beql{lrch}
L_{\La'}(z)\check{R}_k(q^{-2d_k}) = \check{R}_k(q^{-2d_k})L_{\La}(z).
\eeq

Hence we can write
\begin{align}
L_{\La'}(z) D_1\dots D_m \Ec_{\La'}
&=  L_{\La'}(z)\check{R}_k(q^{-2d_k})\check{R}_k(q^{-2d_k})^{-1}D_1\dots D_m\Ec_{\La'}
\non\\
&= \check{R}_k(q^{-2d_k})L_{\La}(z)D_1\dots D_m \check{R}_k(q^{-2d_k})^{-1}\Ec_{\La'} .
\label{prodrld}
\end{align}
The Lemma \ref{lem:et} implies
\beql{eq:RE}
\check{R}_k(q^{-2d_{k}})\Ec_{\La}=\Ec_{\La'}\check{R}_k(q^{2d_{k}}).
\eeq

Hence we can write
\begin{align}
L_{\La'}(z) D_1\dots D_m \Ec_{\La'}
&=  L_{\La'}(z)\check{R}_k(q^{-2d_k})\check{R}_k(q^{-2d_k})^{-1}D_1\dots D_m\Ec_{\La'}
\non\\
&= \check{R}_k(q^{-2d_k})L_{\La}(z)D_1\dots D_m \Ec_{\La}\check{R}_k(q^{-2d_k})^{-1} .
\label{prodrld}
\end{align}

Since the coefficient of $E_{ij}\otimes E_{kl}$ in the $R$-matrix $\check{R}(x)$ is nonzero only if
$\bar i+\bar j +\bar k+\bar l=0\ mod\ 2$,
\ben
\bal
\str^{}_{1,\dots,m}\ts L_{\La'}(z) D_1\dots D_m \Ec_{\La'}
&=\str^{}_{1,\dots,m}\ts\check{R}_k(q^{-2d_k})L_{\La}(z)D_1\dots D_m \Ec_{\La}\check{R}_k(q^{-2d_k})^{-1}\\
&=\str^{}_{1,\dots,m}\ts L_{\La}(z)D_1\dots D_m \Ec_{\La}\check{R}_k(q^{-2d_k})^{-1}\check{R}_k(q^{-2d_k})\\
&=\str^{}_{1,\dots,m}\ts L_{\La}(z) D_1\dots D_m \Ec_{\La}
\eal
\een
thus proving \eqref{lalapr} and completing the proof of Theorem~\ref{thm:centr}.
\epf

Since the series $S_{\La}(z)$ does not depend on the standard $\la$-tableau $\La$,
it is unambiguous to set $S_{\la}(z)=S_{\La}(z)$. We will also use the
formula
\beql{sla}
S_{\la}(z)=\str^{}_{1,\dots,m}\ts T_{\la}(z),\qquad
T_{\la}(z)=\frac{1}{f_{\la}}\sum_{\sh(\La)=\la}\ts L_{\La}(z)D_1\dots D_m\Ec_{\La}.
\eeq

We will now apply Theorem~\ref{thm:centr} to describe
a family of invariants of the
{\em vacuum module at the critical level}
$V_q(\gl_{M|N})$. It is defined as the quotient of $\U_q(\wh\gl_{M|N})_{\text{\rm cri}}$ by the left ideal
generated by all elements $l^-_{ij}[r]$ with $r>0$ and by
the elements $l^-_{i\tss j}[0]-\de_{i\tss j}$ with $i\geqslant j$.
The module $V_q(\gl_{M|N})$ is generated by the vector $\vac$
(the image of $1\in \U_q(\wh\gl_{M|N})_{\text{\rm cri}}$ in the quotient)
such that
$
L^-(u)\tss\vac=I\ts\vac,
$
where $I$ denotes the identity matrix. As a vector space,
$V_q(\gl_{M|N})$ can be identified with the subalgebra
$\Y_q(\gl_{M|N})$ of $\U_q(\wh\gl_{M|N})_{\text{\rm cri}}$ generated by
the coefficients of all series $l^+_{ij}(u)$
subject to the additional relations $l^+_{ii}[0]=1$.
The subspace of invariants of $V_q(\gl_{M|N})$ is defined by
\ben
\z_q(\wh\gl_{M|N})=\{v\in V_q(\gl_{M|N})\ |\ L^-(u)\tss v=I\ts v\}.
\een
One can regard $\z_q(\wh\gl_{M|N})$ as a subspace of $\Y_q(\gl_{M|N})$.
Moreover, this subspace is closed under the multiplication in
the quantum affine algebra. Therefore, $\z_q(\wh\gl_{M|N})$ can be identified
with a subalgebra of $\Y_q(\gl_{M|N})$.

For any standard $\la$-tableau $\La$ introduce the series $\overline S_{\La}(z)$
with coefficients in $\Y_q(\gl_{M|N})$ by
\ben
\overline S_{\La}(z)=\tr^{}_{1,\dots,m}\ts L^+_1(zq^{-2c_1(\La)})\dots L^+_m(zq^{-2c_{m}(\La)})
D_1\dots D_m \Ec_{\La}.
\een

\bco\label{cor:inv}
The series
$\overline S_{\La}(z)\vac$ does not depend on the standard $\la$-tableau $\La$
and all its coefficients
belong to the algebra of invariants $\z_q(\wh\gl_{M|N})$.
Moreover, the coefficients of all series $\overline S_{\La}(z)$ pairwise commute.
\qed
\eco

\bpf It follows from
 Theorem \ref{thm:centr},
that \( L^{-}(u) S_{\La}(z) = S_{\La}(z) L^{-}(u) \).
Applying both sides to the vector
\( 1 \in V_q(\mathfrak{gl}_{M|N}) \)
and using \( S_{\La}(z) 1 = \overline S_{\La}(z) 1 \) we obtain the first part of the corollary.

Applying both sides of the identity \(  S_{\La}(z)   S_{\La'} (w) =S_{\La'} (w)   S_{\La}(z)   \) to the vector \( 1 \),
we have that
\[
 S_{\La}(z)   S_{\La'}(w) 1 =  S_{\La}(z)   \overline S_{\La'}(w) 1
 =
 \overline S_{\La'}(w)S_{\La}(z)1
 = \overline S_{\La'}(w)\overline S_{\La}(z)1
 .
\]
The same calculation for the right hand side gives \( \overline S_{\La'}(w)\overline S_{\La}(z)
=\overline S_{\La}(z)\overline S_{\La'}(w)
\).

\epf

\section{Harish-Chandra images}

To define analogues of the Harish-Chandra homomorphism, we give a total ordering on the generators of $\U_q(\wh{\gl}_{M|N})$.

First, each generator ${l^+_{ij}}^{(r)}$ precedes each generator ${l^-_{kt}}^{(s)}$.
Furthermore, ${l^+_{ij}}^{(r)}\prec {l^+_{kt}}^{(s)}$ if and only if the triple $(i-j,i,r)$
precedes $(k-t,k,s)$ in the lexicographical order. Finally, we set
${l^-_{ij}}^{(r)}\prec {l^-_{kt}}^{(s)}$ if and only if the triple $(j-i,i,r)$
precedes $(t-k,k,s)$ in the lexicographical order. Note that by the RLL relation,
\[
{l^{\pm}_{ij}}^{(r)}{l^{\pm}_{ij}}^{(s)}=(-1)^{\bar i+\bar j}{l^{\pm}_{ij}}^{(s)}{l^{\pm}_{ij}}^{(r)}
\]
for all $r$ and $s$.
Hence, the ordering $\prec$ induces a well-defined total ordering
on the series \eqref{l series} such that $l^+_{ij}(u)\prec l^-_{kt}(u)$ and
\ben
\bal
&l^+_{1,M+N}(u)\prec l^+_{1,M+N-1}(u)\prec\ldots\prec
l^+_{1\tss 1}(u)\prec\ldots\prec
l^+_{M+N,M+N}(u)\prec l^+_{2\tss 1}(u)\prec\ldots\prec l^+_{M+N, 1}(u),
\\[0.3em]
&l^-_{M+N, 1}(u)
\prec l^-_{M+N-1, 1}(u)
\prec   \ldots \prec
l^-_{1\tss 1}(u)\prec\ldots\prec
l^-_{M+N,M+N}(u)\prec l^-_{1\tss 2}(u)\prec\ldots\prec l^-_{1,M+N}(u).
\eal
\een

For the quantum affine superalgebra $\U_q(\wh{\gl}_{M|N})$ at level $c=0$, a PBW basis was first constructed for the standard parity sequence \cite{LWZ} and later generalized to an arbitrary parity sequence in \cite{LZ} for a different sequence.

Consider the ordered monomials in the generators $l_{ij}^{\pm(r)}$. The RLL relations imply that

\beql{relation l(0)}
\bal
 l_{ij}^{\pm}(z)l_{kk}^{+ (0)} =q_i^{-\delta_{ki}}q_j^{\delta_{kj}}l_{kk}^{+(0)}l_{ij}^\pm(z),\\
 l_{ij}^{\pm}(z)l_{kk}^{- (0)} =q_i^{\delta_{ki}}q_j^{-\delta_{kj}}l_{kk}^{-(0)}l_{ij}^\pm(z).
 \eal
 \eeq
Hence we suppose that each monomial only contains either a nonnegative power of $l_{ii}^{+(0)}$ or a positive power of $l_{ii}^{-(0)}$.
Using the same argument given in  \cite{LWZ} and \cite{LZ},
it can be proved that the ordered monomials in the generators form a basis of quantum affine algebra at the critical level $\U_q(\wh\gl_{M|N})_{\text{\rm cri}}$.

We denote by  $L'(z)$  the inverse of $L^{-}(z)$ and its entries by $l'_{ij}(z)$. Then $\U_q(\wh{\gl}_{M|N})$ can be generated by ${l^+_{ij}}^{(r)}$, ${l'_{ij}}^{(r)}$ and ${l^-_{ii}}[0]= {l^+_{ii}}[0]^{-1}$.
Consider the subalgebra  generated by $l^{-(r)}_{ij} , 1\leq i,j\leq M+N$ and $l^{+(0)}_{kk}, 1\leq k\leq M+N $. The mapping
$$L^-(z)\mapsto  L'(z)$$ define an antiautomorphism which maps $ l^{+(0)}_{kk}$ to $ l^{-(0)}_{kk}$.
Then we have a well-defined total ordering on the series such that $l^+_{ij}(u)\prec l'_{kt}(u)$ and
\ben
\bal
&l^+_{1,M+N}(u)\prec l^+_{1,M+N-1}(u)\prec\ldots\prec
l^+_{1\tss 1}(u)\prec\ldots\prec
l^+_{M+N,M+N}(u)\prec l^+_{2\tss 1}(u)\prec\ldots\prec l^+_{M+N, 1}(u),
\\
&{l'}_{M+N, 1}(u)
\prec {l'}_{M+N-1, 1}(u)
\prec   \ldots \prec
{l'}_{1\tss 1}(u)\prec\ldots\prec
{l'}_{M+N,M+N}(u)\prec {l'}_{1\tss 2}(u)\prec\ldots\prec {l'}_{1,M+N}(u).
\eal
\een

Therefore ordered monomials in the generators ${l^+_{ij}}^{(r)}$, ${l'_{ij}}^{(r)}$ that only contains either a nonnegative power of $l_{ii}^{+(0)}$ or a positive power of $l_{ii}^{-(0)}$ form a basis of quantum affine algebra at the critical level $\U_q(\wh\gl_{M|N})_{\text{\rm cri}}$.

Denote by $\U^0$
the subspace of the algebra spanned by
those monomials which do not
contain any generators $l^{+}_{ij}[r]$ and $l'_{ij}[r]$ with $i\ne j$.
Let $x_0$ denote the component
of the linear combination representing the element $x$,
which belongs to $\U^0$.
The mapping $\theta:x\mapsto x_0$ defines the projection
$\theta:\U_q(\wh\gl_{M|N})_{\text{\rm cri}}\to \U^0$. Extending it
by continuity we get the projection $\theta:\wt\U_q(\wh\gl_{M|N})_{\text{\rm cri}}\to \wt\U^0$
to the corresponding completed vector space $\wt\U^0$.

The algebra $\Pi_q(M|N)$ is defined as
the quotient of the algebra of polynomials in independent variables
${l^+_i}^{(r)}$, ${l^-_i}^{(0)}$ and ${l_i'^{(r)}}$ with $i=1,\dots,M+N$ and $r=0,1,\dots$
by the relations ${l^+_i}^{(0)}={l'_i}^{(0)}$, ${l^+_{i}}^{(0)}{l^-_{i}}^{(0)}=1$ for all $i$.
The mapping $\eta:\U^0\to \Pi_q(M|N)$ which takes each ordered monomial
in the generators ${l^{+}_{ii}}^{(r)}$ to the corresponding monomial
in the variables ${l^{+}_{i}}^{(r)}$ , ${l'_{ii}}^{(r)}$ to ${l'_{i}}^{(r)}$, ${l^{-}_{ii}}^{(0)}$ to   ${l^{-}_{i}}^{(0)}$   extend  to an isomorphism
of vector spaces.
Define
the completion $\wt\Pi_q(M|N)$ of the algebra $\Pi_q(M|N)$ as
the inverse limit
\ben
\wt\Pi_q(M|N)=\lim_{\longleftarrow} \Pi_q(M|N)/I_p, \qquad p>0,
\een
where $I_p$ denotes the ideal of $\Pi_q(M|N)$ generated by all elements
$l'_{i}[r]$ with $r\geqslant p$; cf. \eqref{compl}. The isomorphism
$\eta$ extends to an isomorphism of the completed
vector spaces $\eta:\wt\U^0\to \wt\Pi_q(M|N)$.
Thus we get a linear map
\beql{chihom}
\chi:\wt\U_q(\wh\gl_{M|N})_{\text{\rm cri}}\to \wt\Pi_q(M|N)
\eeq
defined as the composition $\chi=\eta\circ\theta$.
The next proposition provides an analogue of the Harish-Chandra homomorphism
for the quantum affine superalgebra.

\bpr\label{prop:hchhom}
The restriction of the map \eqref{chihom} to the center $\Zr_q(\wh\gl_{M|N})$ of the algebra
$\wt\U_q(\wh\gl_{M|N})_{\text{\rm cri}}$ is a homomorphism of commutative algebras
$
\chi:\Zr_q(\wh\gl_{M|N})\to \wt\Pi_q(M|N).
$
\epr

\begin{proof} For \( x, y \in \Zr_q(\wh\gl_{M|N}) \) set \( x_0 = \chi(x) \) and \( y_0 = \chi(y) \). Write \( y \) as a (possibly infinite) linear combination of ordered monomials in the generators \( l_{ij}^{\pm}[r] \). Suppose that
\[
m = \prod_a l_{i_a j_a}^{+(r_a)} \prod_b {l_{i_b j_b}'}^{(r_b)}
\]
is an ordered monomial which occurs in the linear combination.
It follows from Equation \eqref{relation l(0)} that
\beql{index condition}
\sum_a (i_a - j_a) + \sum_b (i_b - j_b) = 0.
\eeq

Since \( x \) is in the center, we have
\[
x m = \prod_a l_{i_a j_a}^{+(r_a)} x \prod_b {l_{i_b j_b}'}^{(r_b)} .
\]
If
$i_a<j_a$ for some $a$ or $i_b<j_b$ for some $b$ then
$\chi(xm)=\chi(m)=0$.

Suppose that $i_a\geq j_a$ and $i_b\geq j_b$  for any $a,b$. It follows from \eqref{index condition} that $i_a= j_a$ and $i_b= j_b$  for any $a,b$.
Then $\chi(m)\neq 0$ if and only if m is of the form

\[
\prod_a l_{i_a i_a}^{+(r_a)}  \prod_b {l_{i_b i_b}'}^{(r_b)}.
\]
Thereofre a nonzero contribution to the image \( \chi(xy) \) can only come from \( \chi(xy_0) \), that is, from expressions of the form

\[
\prod_a l_{i_a i_a}^{+(r_a)} x \prod_b {l_{i_b i_b}'}^{(r_b)}.
\]

Suppose that
\[
p =
\prod_c l_{i_c j_c}^{+(r_c)}  \prod_d {l_{i_d j_d}'}^{(r_d)}
\]
is an ordered monomial which occurs in the linear combination representing \( x \) and \( \chi(p) = 0 \), then $i_c<j_c$ for some $c$ or $i_d<j_d$ for some $d$.

The RTT relation \eqref{RLL} can be written as
\beql{rll2}
\begin{split}
&\delta_{a=c}(zq_a-wq_a^{-1})l_{ab}^{+}(z)l_{cd}^{+}(w)(-1)^{(\bar{a}+\bar{b})(\bar{c}+\bar{d})}+\delta_{a\neq c}(z-w)l_{ab}^{+}(z)l_{cd}^{+}(w)(-1)^{(\bar{a}+\bar{b})(\bar{c}+\bar{d})}\\
&+\delta_{a>c}w(q-q^{-1})l_{cb}^{+}(z)l_{ad}^{+}(w)(-1)^{\bar{c}\bar{b}+\bar{c}\bar{d}+\bar{b}\bar{d}}
+\delta_{a<c}z(q-q^{-1})l_{cb}^{+}(z)l_{ad}^{+}(w)(-1)^{\bar{c}\bar{b}+\bar{c}\bar{d}+\bar{b}\bar{d}}\\
=&\delta_{b=d}(zq_b-wq_b^{-1})l_{cd}^{+}(w)l_{ab}^{+}(z)+\delta_{b\neq d}(z-w)l_{cd}^{+}(w)l_{ab}^{+}(z)\\
&+\delta_{d>b}w(q-q^{-1})l_{cb}^{+}(w)l_{ad}^{+}(z)(-1)^{\bar{c}\bar{b}+\bar{c}\bar{d}+\bar{b}\bar{d}}
+\delta_{d<b}z(q-q^{-1})l_{cb}^{+}(w)l_{ad}^{+}(z)(-1)^{\bar{c}\bar{b}+\bar{c}\bar{d}+\bar{b}\bar{d}}.
\end{split}
\eeq

Suppose that $i_c<j_c$, by the relation \eqref{rll2}, $l_{i i}^{+(r) } l_{i_c j_c}^{+ (r_c)}$ is sum of monomials of the form $\prod_{k=1}^s l_{i_k j_k}^{+(r_k)}$ such that $i_1<j_1$.

For   $i_d<j_d$ , using the relation $R(z/w)^{-1} L'_{1} (z) L'_{2} (w) =   L'_{2} (w)L'_{1} (z) R(z/w)^{-1}$, we can write ${l_{i_d j_d}'}^{(r_c)} {l_{jj}'}^{(r)}$  as sum of monomials of the form $\prod_{k=1}^s {l_{i_k j_k}'}^{(r_k)}$ such that $i_s<j_s$.

Then we conclude that

\[
\chi : \prod_a l_{i_a i_a}^{+(r_a)} p \prod_b {l'_{i_b i_b}}^{(r_b)} \mapsto 0.
\]

Finally, observe that by the defining relations \eqref{rll2}, any two generators \( {l'_{ii}}^{(r)} \) and \( {l'_{jj}}^{(s)}  \) (resp., \( l_{ii}^{+(r)} \) and \( l_{jj}^{+(s)} \)) can be permuted modulo \( \ker \chi \) within any monomial of the form

\[
\prod_a l_{i_a i_a}^{+(r_a)}  \prod_b {l'_{i_b i_b}} ^{  (r_b)} .
\]

This proves that \( \chi(xy) = x_0 y_0 \).

\end{proof}

Combine the generators of the algebra $\Pi_q(M|N)$ into the series
\begin{equation}\label{l series}
\begin{split}
l^{+}_{i}(u)= \sum_{r=0}^{\infty} {l^{+}_{i}}^{(r)}u^{ r}, \quad l'_{i}(u)= \sum_{r=0}^{\infty} {l'_{i}}^{(r)}u^{ -r},
\end{split}
\end{equation}
and for $i=1,\dots,M+N$ set
\begin{align}
x_i(z)=\left\{ \begin{aligned}
& q^{2i} l^{+}_{i}(z)\ts l'_{i}(q^{M-N}z)    ,\ & 1\leq i\leq M,\\
&q^{2M+2-2i} l^{+}_{i}(z)\ts l'_{i}(q^{M-N}z), \ & M+1\leq i\leq M+N.\\
\end{aligned} \right.
\end{align}

This is a Laurent series in $z$ whose coefficients
are elements of the completed algebra $\wt\Pi_q(M|N)$.
We are now in a position to state the main result of this section.

\bth\label{thm:hchim}
The image of the series $S_{\la}(z)$ under the Harish-Chandra homomorphism is as follows
\ben
\chi:S_{\la}(z)\mapsto \sum_{\sh(\mathcal{T})=\la} (-1)^{\overline{\mathcal{T}}}
\prod_{\al\in \la}x^{}_{\mathcal{T}(\al)}(zq^{-2c(\al)}),
\een
summed over semistandard tableau $\mathcal{T}$ of shape $\la$ with entries in $\{1,2,\dots, M+N\}$，where $\overline{\mathcal{T}}
=\sum_{\alpha\in\lambda}
\overline{x_{\mathcal{T}(\al)}}$.
\eth

The proof of Theorem~\ref{thm:hchim} will be given
in the rest of this section.

For any $m$-tuple $I=(i_1,\dots, i_m )$ with
$1\leqslant i_1\leqslant \dots \leqslant i_m\leqslant M+N$,
we simply denote  $I=(1^{m_1},\ldots,(M+N)^{m_{M+N}})$, where $m_i=m_i(I)= \mathrm{Card}\{j\in I |j=i\}$.
We thus get the corresponding composition
$\mu_{I}=(m_1,\dots, m_{M+N})$ of $m$.
Let $\Sym_{I}$ be the  Young subgroup $ \Sym_{m_1}\times \dots \times \Sym_{m_{M+N}}$ of
the symmetric group $\Sym_m$ and $\mathcal M_{I}$ be the minimal length coset representative of $\Sym_m/\Sym_I$.
Given $\si\in\Sym_m$, the unique decomposition $\si=\om\tss\pi$
with $\om\in\mathcal M_{I}$ and $\pi\in\Sym_{I}$.
Moreover,
by counting the number of inversions in the permutation $\si$ we find that
$l(\si)=l(\om)+l(\pi)$. This implies the relation $T_{\si}=T_{\om}T_{\pi}$
in the Hecke algebra $\Hc_m$.

For any ring $R$, it will be convenient to use a standard notation for the matrix elements
$A^{i_1\dots i_m}_{j_1\dots j_m}$ of an even operator
\ben
A=\sum_{i_1,\dots,i_m,\ts j_1,\dots,j_m}A^{i_1\dots i_m}_{j_1\dots j_m} \ot e_{i_1,j_1}\ot \dots \ot e_{i_m,j_m}\in R\ot (\End\CC^{M|N})^{\ot m},
\een
where $e_{ij}$ is the standard basis of $\mathrm{End}(\mathbb{C}^{M|N})$. Then we can write the bases of $(\mathbb C^{M|N})^{\ot m}$  and its dual bases of $((\mathbb C^{M|N})^{\ot m})^{*}$  respectively:
\beq
\mid i_1,\dots,i_m\rangle =e_{i_1}\otimes\cdots \otimes e_{i_m}, \qquad \langle i_1,\dots ,i_m\mid=(e_{i_1}\otimes\cdots \otimes  e_{i_m})^*,
\eeq
with parity $\bar I=\bar i_1+\ldots +\bar i_m$.
Then the  coefficients of the even  operator $A$ are given by
\beq
A_{j_1,\ldots,j_m}^{i_1,\ldots,i_m}=(-1)^{\gamma(I,J)}\langle i_1,\ldots ,i_m\mid A \mid j_1,\ldots,j_m\rangle,
\eeq
where
$\gamma(I,J)=\sum_{a} \bar i_a(\bar j_a+1)+\sum_{a<b}  \bar j_b(\bar i_a +\bar j_a)$. In particular,
\beq
\bal
A_{i_1,\ldots,i_m}^{i_1,\ldots,i_m}&=\langle i_1,\ldots ,i_m\mid A \mid i_1,\ldots,i_m\rangle,\\
str_{1,\ldots,m}(A)&=\sum_{I}(-1)^{\bar I}\langle i_1,\ldots ,i_m\mid A \mid i_1,\ldots,i_m\rangle.
\eal
\eeq

Recall that
\begin{eqnarray}
\begin{array}{rcl}
R &=&  \sum\limits_{i=1}^{M+N} q_i   E_{ii} \otimes E_{ii}  +   \sum\limits_{i \neq j} E_{ii} \otimes E_{jj} +  \sum\limits_{i<j}(q_j-q_j^{-1})  E_{ij} \otimes E_{ji}.
\end{array}
\end{eqnarray}

Setting $\check{R}=PR$, we get
\ben
\check{R}_{k}\check{R}_{k+1}\check{R}_{k}=\check{R}_{k+1}\check{R}_{k}\check{R}_{k+1}\Fand
(\check{R}_{k}-q)(\check{R}_{k}+q^{-1})=0,
\een
where $\check{R}_{k}=P_{k,k+1}R_{k,k+1}$.
We have the property $\check{R}(e_i\ot e_j)=(-1)^{\bar i\bar j}e_j\ot e_i$
for $i<j$.
For an $m$-tuple $I=(i_1,\dots, i_m)$ with
$1\leqslant i_1\leqslant \dots \leqslant i_m\leqslant M+N$ and $\om\in \mathcal M_{I}$,
set
\ben
\overline{I_{\om}}=\sum_{j<k, \atop \om(j)> \om(k)}\bar i_j \bar i_k.
\een

Then we have
\ben
\check{R}_{\om}\tss|\tss i_1,\dots, i_m\rangle=(-1)^{\overline{I_{\om}}}|\tss i_{\om(1)},\dots, i_{\om(m)}\rangle
\een
and
\ben
\langle i_1,\dots, i_m\tss |\tss\check{R}_{\om^{-1}}=(-1)^{\overline{I_{\om}}}\langle i_{\om(1)},\dots, i_{\om(m)}\tss|.
\een
Formula \eqref{sla} for $S_{\la}(z)$ then implies
\begin{multline}
S_{\la}(z)
=\sum_{j_1,\dots, j_m=1}^{M+N}(-1)^{\bar J}\langle j_1,\dots, j_m\tss|\tss T_{\la}(z)\tss |\tss j_1,\dots,j_m\rangle\\
=\sum_{i_1\leqslant\dots\leqslant i_m} (-1)^{\bar I}
\sum_{\om\in \mathcal M_{I}}\langle i_1,\dots, i_m\tss|\tss
\check{R}_{\om^{-1}}T_{\la}(z)\check{R}_{\om} \tss|\tss i_1,\dots, i_m\rangle.
\non
\end{multline}

Using Lemmas\ref{lem:et}, \ref{lem:LE} and \eqref{isomhecke}, we can prove the following lemma.
\ble\label{lem:RT}
For any $\si\in \Sym_m$ and $\la \vdash m$ we have
\ben
\check{R}_{\si}{T}_{\la}(z)={T}_{\la}(z)\check{R}_{\si}.
\een
\ele

Lemma~\ref{lem:RT} allows us to write the expression for
$S_{\la}(z)$ in the form
\ben
\bal
S_{\la}(z)&=\sum_{i_1\leqslant\dots\leqslant i_m} (-1)^{\bar I}
\sum_{\om\in \M_{I}}\langle i_1,\dots, i_m\tss|\tss T_{\la}(z)
\check{R}_{\om^{-1}}\check{R}_{\om} \tss|\tss i_1,\dots, i_m\rangle\\
&=\frac{1}{f_{\la}}\sum_{i_1\leqslant\dots\leqslant i_m}
(-1)^{\bar I}
\sum_{\om\in \M_{I}}\sum_{\sh(\La)=\la}
\langle i_1,\dots, i_m\tss|\tss L_{\La}(z)D_1\dots D_m
\Ec_{\La}\check{R}_{\om^{-1}}\check{R}_{\om}\tss|\tss i_1,\dots, i_m\rangle.
\eal
\een

By \eqref{tauchi},
we have that
\ben
\Ec_{\La}\check{R}_{\om^{-1}}\check{R}_{\om}
=\sum\limits_{\si\in \Sym_m}\tau\left(e_{\La}T_{\om^{-1}}T_{\om}T_{\si^{-1}}\right)\check{R}_{\si}
=\sum\limits_{\si\in \Sym_m}\frac{1}{c_{\la}}
\chi_{\la}(e_{\La}T_{\om^{-1}}T_{\om}T_{\si^{-1}})\check{R}_{\si}.
\een
Therefore, $S_{\la}(z)$ equals
\ben
\frac{1}{f_{\la}{c_{\la}}}\sum_{i_1\leqslant\dots\leqslant i_m}
(-1)^{\bar I}
\sum_{\substack{\om\in \M_{I}\\\si\in \Sym_m}}\sum_{\sh(\La)=\la}
\chi_{\la}\left(e_{\La}T_{\om^{-1}}T_{\om}T_{\si^{-1}}\right)
\langle i_1,\dots, i_m\tss|\tss L_{\La}(z)D_1\dots D_m
\check{R}_{\si}\tss|\tss i_1,\dots, i_m\rangle.
\een
Using the decomposition $\Sym_m=\M_{I}\Sym_{I}$, we can write this expression
in the form
\begin{multline}
\frac{1}{f_{\la}{c_{\la}}}\sum_{i_1\leqslant\dots\leqslant i_m}(-1)^{\bar I}
\sum_{\substack{\om,\om'\in \M_{I}\\\pi\in \Sym_{I}}}\sum_{\sh(\La)=\la}
\chi_{\la}\big(e_{\La}T_{\om^{-1}}T_{\om}T_{\pi^{-1}}T_{\om'^{-1}}\big)\\[-1em]
{}\times\langle i_1,\dots, i_m\tss|\tss L_{\La}(z)D_1\dots D_m
\check{R}_{\om'}\check{R}_{\pi}\tss|\tss i_{1},\dots, i_{m}\rangle.
\non
\end{multline}
We denote $\Sym_{0}=\Sym_{ m_1} \times \Sym_{ m_2} \times \cdots \times \Sym_{ m_M}$,
and  $\Sym_{1}=\Sym_{ m_{M+1}} \times \Sym_{m_{M+2}} \times \cdots \times \Sym_{ m_{M+N}}$.
Then $\Sym_{I}=\Sym_{0}\times \Sym_{1}$.
For an $m$-tuple $(i_1,\dots, i_m)$ with $1\leqslant i_1\leqslant \dots \leqslant i_m \leqslant M+N$, $\pi\in \Sym_{I}$
can be written as $\pi_0\pi_1$ where $\pi_i\in\Sym_{i}$. One easily verifies that
\ben
\check{R}_{\pi}\tss|\tss i_{1},\dots, i_{m}\rangle
=q^{l(\pi_0)}(-q)^{-l( \pi_1)}\tss|\tss i_{1},\dots, i_{m}\rangle.
\een
Moreover, we have
\ben
\bal
\check{R}_{\om'}\check{R}_{\pi}\tss|\tss i_{1},\dots, i_{m}\rangle
=q^{l(\pi_0)}(-q)^{-l( \pi_1)}\check{R}_{\om'}\tss|\tss i_{1},\dots, i_{m}\rangle\\
=(-1)^{ \overline{I_{\om'}}} q^{l(\pi_0)}(-q)^{-l( \pi_1)}\tss|\tss i_{\om'(1)},\dots, i_{\om'(m)}\rangle.
\eal
\een
Hence
\begin{multline}
S_{\la}(z)=\frac{1}{f_{\la}{c_{\la}}}\sum_{i_1\leqslant\dots\leqslant i_m}
(-1)^{\bar I}
\sum_{\substack{\om,\om'\in \M_{I}\\\pi\in \Sym_{I}}}\sum_{\sh(\La)=\la}
 \chi_{\la}\big(e_{\La}T_{\om^{-1}}T_{\om}T_{\pi^{-1}}T_{\om'^{-1}}\big)\\
{}\times (-1)^{\overline{I_{\om'}}} q^{l(\pi_0)}(-q)^{-l( \pi_1)}\langle i_1,\dots, i_m\tss|\tss L_{\La}(z)D_1\dots D_m
\tss|\tss i_{\om'(1)},\dots, i_{\om'(m)}\rangle
\non
\end{multline}
which equals
\begin{multline}
\frac{1}{f_{\la}c_{\la}}\sum_{i_1\leqslant\dots\leqslant i_m}(-1)^{\bar I}
\sum_{\substack{\om,\om'\in \M_{I}\\\pi\in \Sym_{I}}}
\sum_{\sh(\La)=\la}\sum_{j_1,\dots, j_m}
q^{l(\pi_0)}(-q)^{-l( \pi_1)}\chi_{\la}\big(e_{\La}
T_{\om^{-1}}T_{\om}T_{\pi^{-1}}T_{\om'^{-1}}\big)\\
{}\times (-1)^d\ l_{i_{1},j_1}^{+}(zq^{2c_{1}(\La)})\dots l_{i_{m},j_m}^{+}(zq^{2c_{m}(\La)})
\ts \wt{l}_{j_m ,i_{\om'(m)}}(zq^{M-N+2c_{m}(\La)})
 \dots \wt{l}_{j_1 ,i_{\om'(1)}}(zq^{M-N+2c_{1}(\La)}),
\non
\end{multline}
where we used the entries of the matrix $L^-(z)^{-1}D=[\tss\wt l_{ij}(z)]$ and $d$ is as follows:
\beq
d=    \overline{I_{\om'}}  +\sum\limits_{a<b}\bar {j_a}(\overline{i_b}+\overline{i_{\om'(b)}}).
\eeq

The definition of the homomorphism \eqref{chihom} implies that
a nonzero contribution to the Harish-Chandra
image of $S_{\la}(z)$ comes only from
 the
summands  with $i_1\geq j_1\geq i_{\om'(1)}$. This condition implies that $\om'(1)=1$ and $i_1=j_1$. The same observation gives $\om'(2)=2$ and $i_2=j_2$,\ldots, $\om'(m)=m$ and $i_m=j_m$.
In these case, $d=0 \ mod\ 2$.
Therefore, the image of $S_{\la}(z)$ is given by
\begin{multline}
\frac{1}{f_{\la}c_{\la}}\sum_{i_1\leqslant\dots\leqslant i_m}
(-1)^{\bar I}
\sum_{\substack{\om\in \M_{I}\\\pi\in \Sym_{I}}}\sum_{\sh(\La)=\la}
q^{l(\pi_0)}(-q)^{-l( \pi_1)}
\chi_{\la}\big(e_{\La}T_{\om^{-1}}T_{\om}T_{\pi^{-1}}\big)\\
{}\times l_{i_1,i_1}^{+}(zq^{2c_{1}(\La)})\dots l_{i_{m},i_m}^{+}(zq^{2c_{m}(\La)})
\ts\wt{l}_{i_m ,i_{m}}(zq^{M-N+2c_{m}(\La)})\dots
 \wt{l}_{i_1 ,i_{1}}(zq^{M-N+2c_{1}(\La)})
\non
\end{multline}
which coincides with
\ben
\frac{1}{f_{\la}c_{\la}}\sum_{i_1\leqslant\dots\leqslant i_m}
(-1)^{\bar I}
\sum_{\substack{\om\in \M_{I}\\\pi\in \Sym_{I}}}\sum_{\sh(\La)=\la}
q^{l(\pi_0)}(-q)^{-l( \pi_1)}  \chi_{\la}\big(e_{\La}T_{\om^{-1}}T_{\om}T_{\pi}\big)
x_{i_1}(zq^{2c_{1}(\La)})\dots x_{i_m}(zq^{2c_{m}(\La)}).
\een

Given an $m$-tuple $(i_1,\dots, i_m)$ with $i_1\leqslant i_2\leqslant \dots \leqslant i_m$,
we let $\Tc=i(\La)$ denote the
tableau obtained from a standard $\la$-tableau $\La$ by replacing
the entry $r$ with $i_r$ for $r=1,\dots, m$.
We denote $\overline{i_1}+\cdots+\overline{i_m}$ by $\overline{\Tc}$.
The entries of $\Tc$ weakly increase
along the rows and down the columns.
Changing the order of summation in the above expression
for the Harish-Chandra image of $S_{\la}(z)$, we can write it as
\ben
\bal
\frac{1}{f_{\la}{c_{\la}}}
\sum_{\Tc,\ts\sh(\Tc)=\la}
(-1)^{\overline{\Tc}}
\Bigg(
\sum_{\substack{\sh(\La)=\la\\ }}\ts
\sum_{\substack{\om\in \M_{I}\\\pi\in \Sym_{I}}}
q^{l(\pi_0)}(-q)^{-l( \pi_1)}
\chi_{\la}\big(e_{\La}T_{\om^{-1}}T_{\om}T_{\pi}\big)
\Bigg)\prod_{\al\in \la}x^{}_{\Tc(\al)}(zq^{2c(\al)}),
\eal
\een
where $\Tc$ runs over $\la$-tableaux with entries in $\{1,\dots,M+N\}$ such that
the entries of $\Tc$ weakly increase
along the rows and down the columns.
To complete the proof of Theorem~\ref{thm:hchim}, it will therefore be enough to show that
for a given $\Tc$ we have
\beql{idenchi}
\sum_{\substack{\sh(\La)=\la\\ }}\ts
\sum_{\substack{\om\in \M_{I}\\\pi\in \Sym_{I}}}
{ q^{l(\pi_0)}(-q)^{-l( \pi_1)}}
\chi_{\la}\big(e_{\La}T_{\om^{-1}}T_{\om}T_{\pi}\big)
=\begin{cases}f_{\la}\tss c_{\la}\quad&\text{if $\Tc$ is semistandard},\\
     0\quad &\text{otherwise.}
   \end{cases}
\eeq
Let $\Hc_{I}$ denote the subalgebra of $\Hc_m$ corresponding to
the Young subgroup $\Sym_{I}$ of $\Sym_m$.
Introduce the notation
\ben
s_{I}=\sum\limits_{\pi\in \Sym_{I}}q^{l(\pi_0)}(-q)^{-l( \pi_1)}T_{\pi}\in\Hc_{I}\Fand
e^{}_{\Tc}=\sum\limits_{\substack{\sh(\La)=\la }\atop i(\Lambda)=\Tc}e^{}_{\La}.
\een
 
The element $s_{I}$ is proportional to an idempotent
\beql{misq}
s_{I}^2=\prod_{r=1}^M{[\al_r]_q!\ts q^{m_r(m_r-1)/2}}\prod_{r=M+1}^{M+N}{[\al_r]_{q^{-1}}!\ts q^{-m_r(m_r-1)/2}}\ts s_{I}.
\eeq

The left hand side of \eqref{idenchi} can now be written as
\beql{chico}
\sum\limits_{\om\in \M_{I}} \chi_{\la}\big(e_{\Tc}T_{\om^{-1}}T_{\om}s_{I}\big)=
\sum\limits_{\om\in \M_{I}} \chi_{\la}\big(T_{\om}s_{I}e_{\Tc}T_{\om^{-1}}\big).
\eeq
We will use the homomorphism $\vp_{\la}$ associated with the $\Hc_m$-module $V_{\la}$.

\ble\label{lem:ASchur}
We have the relation
\beql{phichi}
\sum_{\om\in \M_I}\vp_{\la}\big(T_{\om}s_{I}e_{\Tc}T_{\om^{-1}}\big)
=\prod_{r=1}^M\frac{1}{[\al_r]_q!\ts q^{\al_r(\al_r-1)/2}}\ts
\prod_{r=M+1}^{M+N}\frac{q^{\al_r(\al_r-1)/2}}{[\al_r]_{q^{-1}}!}\ts
c_{\la}\ts\chi_{\la}(s_{I}e_{\Tc})\ts{\rm id}^{}_{V_{\la}}.
\eeq
\ele

\bpf
Denote by
$V_{\Tc}$ the linear span of the basis vectors $v_{\La}\in V_{\la}$ such that
$i(\La)=\Tc$. It follows from \eqref{eq:Young basis} that the subspace $V_{\Tc}$
of $V_{\la}$ is invariant under
the action of the subalgebra $\Hc_{I}$.
Therefore, if a standard $\la$-tableau $\La$ is such that $i(\La)=\Tc$, then
$\vp_{\la}(e_{\Tc}s_{I})(v_{\La})=\vp_{\la}(s_{I})(v_{\La})$.
Otherwise, if $i(\La)\neq \Tc$ then
\ben
\bal
\vp_{\la}(e_{\Tc}s_{I})(v_{\La})&=\sum_{\sh(\Ga)=\la}
\langle s_{I}v_{\La},v^{}_{\Ga}\rangle e^{}_{\Tc}v^{}_{\Ga}
=\sum_{\sh(\Ga)=\la}\langle v_{\La},s_{I}^{*}v^{}_{\Ga}\rangle
e^{}_{\Tc}v^{}_{\Ga}\\
&=\sum\limits_{\substack{\sh(\Ga)=\la\\i(\Ga)=\Tc}}
\langle v_{\La},s_{I}v^{}_{\Ga}\rangle e^{}_{\Tc}v^{}_{\Ga}
+\sum\limits_{\substack{\sh(\Ga)=\la\\i(\Ga)\ne\Tc}}
\langle v_{\La},s_{I}v^{}_{\Ga}\rangle e^{}_{\Tc}v^{}_{\Ga}=0.
\eal
\een
On the other hand, we also have
\ben
\vp_{\la}(s_{I}e_{\Tc})(v_{\La})=\begin{cases}
   0\qquad&\text{if\ \  $i(\La)\neq \Tc$},\\
    \vp_{\la}( s_{I})(v_{\La})\qquad&\text{if\ \  $i(\La)= \Tc$.}
   \end{cases}
\een
Therefore, we may conclude that $\vp_{\la}(s_{I}e_{\Tc})=\vp_{\la}(e_{\Tc}s_{I})$.

Furthermore, taking into account \eqref{misq}, we can write
\ben
\sum_{\om\in \M_I}\vp_{\la}(T_{\om}s_Ie_{\Tc}T_{\om^{-1}})
=\prod_{r=1}^M\frac{1}{[\al_r]_q!\ts q^{\al_r(\al_r-1)/2}}\prod_{r=M+1}^{M+N}\frac{1}{[\al_r]_{q^{-1}}!\ts q^{-\al_r(\al_r-1)/2}}\ts
\sum_{\om\in \M_I}\vp_{\la}\big(T_{\om}s_Ie_{\Tc}s_IT_{\om^{-1}}\big).
\een
Since $T_{\pi}s_I=q^{l(\pi_0)}(-q)^{-l( \pi_1)}s_I$ for $\pi\in \Sym_{I}$, we have
\begin{multline}
\sum_{\om\in \M_I}\vp_{\la}\big(T_{\om}s_Ie_{\Tc}s_IT_{\om^{-1}}\big)
=\sum\limits_{\substack{\om\in \M_I\\\pi\in \Sym_{I}}}
\vp_{\la}\big( q^{l(\pi_0)}(-q)^{-l(\pi_1)} T_{\om}s_Ie_{\Tc}T_{\pi^{-1}}T_{\om^{-1}}\big)\\[-1em]
{}=\sum\limits_{\substack{\om\in \M_I\\\pi\in \Sym_{I}}}
\vp_{\la}\big(T_{\om}T_{\pi}s_Ie_{\Tc}T_{\pi^{-1}}T_{\om^{-1}}\big).
\non
\end{multline}
Now observe that
$\{T_{\om}T_{\pi}\tss|\tss\om\in \M_I, \pi\in \Sym_{I}\}$
and $\{T_{\pi^{-1}}T_{\om^{-1}}\tss|\tss\om\in \M_I, \pi\in \Sym_{I}\}$
are dual bases of $\Hc_m$.
Hence by the fact \eqref{iuschur}, we can get
\ben
\sum_{\om\in \M_I}\vp_{\la}\big(T_{\om}s_Ie_{\Tc}T_{\om^{-1}}\big)
=\prod_{r=1}^M\frac{1}{[\al_r]_q!\ts q^{\al_r(\al_r-1)/2}}\prod_{r=M+1}^{M+N}\frac{1}{[\al_r]_{q^{-1}}!\ts q^{-\al_r(\al_r-1)/2}}\ts
c_{\la}\ts\tr^{}_{V_{\la}}\big(\vp_{\la}(s_Ie_{\Tc})\big)
{\rm id}^{}_{V_{\la}},
\een
which equals the right hand side of \eqref{phichi}.
\epf

Returning to the proof of \eqref{idenchi}, note that by Lemma~\ref{lem:ASchur} the expression
\eqref{chico} equals
\ben
\sum\limits_{\om\in \M_I} \chi_{\la}
\left(T_{\om}s_Ie_{\Tc}T_{\om^{-1}}\right)
=\prod_{r=1}^M\frac{1}{[\al_r]_q!\ts q^{\al_r(\al_r-1)/2}}\prod_{r=M+1}^{M+N}\frac{1}{[\al_r]_{q^{-1}}!\ts q^{-\al_r(\al_r-1)/2}}\ts
 c_{\la}\tss f_{\la}\ts
\chi_{\la}(s_Ie_{\Tc}).
\een
Thus, the following lemma will complete the proof of Theorem~\ref{thm:hchim}.

\ble\label{lem:skew}
We have
\ben
\chi_{\la}(s_Ie_{\Tc})
=\begin{cases} \prod\limits_{r=1}^M [\al_r]_{q }!\ts q^{\al_r(\al_r-1)/2}\prod\limits_{r=M+1}^{M+N} [\al_r]_{q^{-1}}!\ts q^{-\al_r(\al_r-1)/2} \ts
\quad&\text{if $\Tc$ is semistandard},\\
     0\quad &\text{otherwise.}
     \end{cases}
\een
\ele

\bpf
Write
\ben
s_Ie_{\Tc}=
\sum\limits_{\substack{\pi\in \Sym_{I}}}{ q^{l(\pi_0)}(-q)^{-l( \pi_1)}}T_{\pi}e_{\Tc}
=
\prod_{r=1}^{M+N}\sum\limits_{\substack{\si_r\in \Sym_{m_r}}} ((-1)^{\bar r}q)^{(-1)^{\bar r}l(\si_r)}     T_{\si_r}e_{\Tc}.
\een
Hence
\ben
\bal
\chi_{\la}(s_Ie_{\Tc})
=\sum\limits_{\substack{\sh(\La)=\la\\ }}
\prod_{r=1}^{M+N}\sum\limits_{\substack{\si_r\in \Sym_{m_r}}}  ((-1)^{\bar r}q)^{(-1)^{\bar r}l(\si_r)}
\langle T_{\si_r}v_{\La},v_{\La}\rangle\\
=\prod _{r=1}^{M+N}\sum_{\si_r\in \Sym_{\al_{r}}} ((-1)^{\bar r}q)^{(-1)^{\bar r}l(\si_r)}
\chi_{\om_r}(T_{\si_{r}}),
\eal
\een
where $\om_r$ is the skew diagram which
consists of the boxes of $\Tc$ occupied by $r$, and
$\chi_{\om_r}$ denotes the skew character of $\Hc_{m_r}$ associated with $\om_r$.
Since $q^{l(\si_r)}=\chi_{\iota}(T_{\si_{r}^{-1}})$(resp. $(-q)^{-l(\si_r)}=\chi_{\varepsilon}(T_{\si_{r}^{-1}})$) for the trivial(resp. sign)
representation $\iota$(resp. $\varepsilon$) of the Hecke algebra $\Hc_{m_r}$,
we can write
\ben
\chi_{\la}(s_Ie_{\Tc})=\prod _{r=1}^{M}\sum_{\si_r\in \Sym_{\al_{r}}}\chi_{\iota}
(T_{\si_{r}^{-1}})
\chi_{\om_r}(T_{\si_{r}}).
\prod _{r=M+1}^{M+N}\sum_{\si_r\in \Sym_{\al_{r}}}\chi_{\varepsilon}
(T_{\si_{r}^{-1}})
\chi_{\om_r}(T_{\si_{r}}).
\een
The multiplicity of the trivial representation
$\iota$ in the skew representation of $\Hc_{\alpha_{r}}$
associated with $\om_{r} $ is zero unless $\om_{r}$
does not contain two boxes in the same column, in which case the multiplicity is $1$.
Then by \eqref{orthcha},
\ben
\sum_{\si_r\in \Sym_{\al_{r}}}q^{l(\si_r)}
\chi_{\om_r}(T_{\si_{r}})
=\sum_{\si_r\in \Sym_{\al_{r}}}
\chi_{\iota}(T_{\si_{r}})\chi_{\iota}(T_{\si_{r}^{-1}})
=[\al_r]_q!\ts q^{\al_r(\al_r-1)/2}.
\een
\ben
\sum_{\si_r\in \Sym_{\al_{r}}}(-q)^{-l(\si_r)}
\chi_{\om_r}(T_{\si_{r}})
=\sum_{\si_r\in \Sym_{\al_{r}}}
\chi_{ \varepsilon}(T_{\si_{r}})\chi_{\iota}(T_{\si_{r}^{-1}})
=[\al_r]_{q^{-1}}!\ts q^{-\al_r(\al_r-1)/2}.
\een
This proves the lemma and completes the proof of Theorem~\ref{thm:hchim}.
\epf

\bigskip
\centerline{\bf Acknowledgments}
\medskip

The work is supported in part by the National Natural Science Foundation of China (grant nos. 12171303, 12471026, and 12571026) and the Simons Foundation (grant no. MP-TSM-00002518).

\newpage
\noindent
N.J.:\newline
Department of Mathematics\\
North Carolina State University, Raleigh, NC 27695, USA\\
jing@ncsu.edu

\vspace{5 mm}

\noindent
M.L.:\newline
School of Artificial Intelligence\\
Jianghan University,\\
 Wuhan 430056, Hubei, China\\
ming.l1984@gmail.com

\vspace{5 mm}

\noindent
J.Z.:\newline
School of Mathematics and Statistics\\
Central China Normal University,\\
Wuhan, Hubei 430079, China\\
jzhang@ccnu.edu.cn

\end{document}